\documentclass[12pt]{article}
\usepackage{pgf,tikz}
\usetikzlibrary{arrows}
\usepackage{amsmath}%
\usepackage{amsfonts}%
\usepackage{amssymb}%
\usepackage{graphicx}%
\newtheorem{theorem}{Theorem}

\newtheorem{lemma}{Lemma}

\begin{document}
\title{On Star critical Ramsey numbers related to  large cycles versus complete graphs }
\author{Chula J. Jayawardene \\
Department of Mathematics\\
University of Colombo \\
Sri Lanka\\
email: c\_jayawardene@maths.cmb.ac.lk\\
\\
W. Chandanie W. Navaratna \\
Department of Mathematics \\
The Open University of Sri Lanka \\
Sri Lanka\\
email: wcper@ou.ac.lk \\
}

\maketitle
\begin{abstract} Let $K_n$ denote the complete graph on $n$ vertices and $G, H$ be finite graphs without loops or multiple edges. Consider a  two-coloring of edges of $K_n$. When a copy of $G$ in the first color, red, or a copy of $H$ in the second color, blue is in $K_n$, we write  $K_n\rightarrow (G,H)$. The Ramsey number $r(G, H)$ is defined as the smallest positive integer $n$ such that $K_{n} \rightarrow (G, H)$.  Star critical Ramsey $r_*(G, H)$ is defined as the largest value of $k$  such that $K_{r(G,H)-1} \sqcup K_{1,k} \rightarrow (G, H)$.
In this paper, we find $r_*(C_n, K_m)$ for $m \geq 7$ and  $n \geq (m-3)(m-1)$.  
\end{abstract}


\section{Introduction}
\noindent Star-critical Ramsey numbers introduced by Hook and Isaak \cite{HoIs,Ho} in 2010 have captured the attention of many authors in the recent years. Literature reveals calculation of Ramsey numbers related to  $r(C_n, K_m)$ for $n \geq m$  and $m \leq 7$ (see \cite{Is,Ra}) and Star-critical Ramsey numbers related to  $r(C_n, K_3)$ for $n \geq 3$, $r(C_n, K_4)$ for $n \geq 4$, $r(C_n, K_5)$ for $n \geq 5$ and $r(C_n, K_6)$ for $n \geq 15$. In this paper, we extend the calculation of Star-critical Ramsey numbers to cover $r(C_n, K_m)$ for $m \geq 7$ and  $n \geq (m-1)(m-3)$. In particular, we show that $r_*(C_n, K_m)=(m-2)(n-1)+2$ for  $m \geq 7$ and $n \geq (m-3)(m-1)$.

\vspace{14pt}

\section{Notation}
\noindent For ease of reference, we borrow the notation used in \cite{Ja,JaRo,JaRo1}. Given a graph $G$, we say $Y \subseteq V(G)$ is an \textit{independent set} if no pair  of  vertices of $Y$ is adjacent to each other in $G$. Equivalently, $Y$ forms a clique in $G^c$. The independence number $\alpha(G)$ is defined as the size of the largest independent set. Thus, $\alpha(G)=\max \{ \left| I\right| : I$ is an independent set of $G \}$. Given a graph $G$ and a vertex $v \in V(G)$, we define the \textit{neighbourhood of $v$ in $G$},  $\Gamma(v)$, as the set of vertices adjacent to $v$ in $G$. The \textit{degree} of a vertex $v$, $d(v)$, is defined as the cardinality of  $ $ $\Gamma(v)$, $ $ i.e. $d(v)=| \Gamma(v)|$. We write $\bar \Gamma(v)$ for $\Gamma(v) \cup  \{ v \}$. The minimum degree of a graph $G(V,E)$, denoted by $\delta(G)$, is defined as  $ $ $ \min\{ d(v) | v\in V  \}$.  Given a graph $G$ and a non-empty subset $S$ of $V$, \textit{the induced subgraph of $S$ in $G$} denoted by $G[S]$ is defined as the  subgraph obtained by deleting all the vertices of $S^c$ from $G$. Moreover, $G \setminus S$ is defined as $G[V(G) \setminus S]$. Given a graph $G$ and two disjoint subgraphs $H$ and $K$ of $G$, we denote the set of edges between $H$ and $K$ by $E(H,K)$ and define $G-H$ as the subgraph of $G$ obtained by deleting all the edges of $H$ from $G$.

\vspace{10pt}

\section{Useful lemmas for calculating $r_*(C_n, K_m)$ for $m \geq 7$ and $n \geq (m-3)(m-1)$}

We first present four lemmas that were used to to calculate $r_*(C_n, K_m)$ when $m \geq 7$ and $n \geq (m-3)(m-1)$ of which Lemma \ref{l1} is from \cite{JaRo}, Lemma  \ref{l2} is from \cite{Ja} and Lemma  \ref{l3} is from \cite{BoJaSh} by Bollab{\'a}s et al. We provide an inductive proof for Lemma \ref{l4}.
\vspace{10pt}

\begin{lemma}
\label{l1}
\noindent (\cite{JaRo}, Lemma 2). A $C_n$- free graph $G$ of order $N$  with independent number  less than or equal to $m$ has minimal degree  greater than or equal to $N-r(C_n,K_{m})$.
\end{lemma}

\vspace{2pt}

\begin{lemma}
\label{l2}
\noindent (\cite{Ja}, Lemma 8). A $C_n$ -free graph (where $n \geq 15$) of order $5(n-1)$ with no independent set of 6 vertices contains a $5K_{n-1}$.
\end{lemma}

\vspace{2pt}

\begin{lemma} (\cite{BoJaSh}, Lemma 5)
\label{l3}

\noindent Suppose $G$ contains the cycle $(u_1, u_2, ... , u_{n-1}, u_1)$ of length $n - 1$ but
no cycle of length $n$. Let $Y = V(G) \setminus \{u_1, u_2, ... , u_{n-1}\}$. Then,

\vspace{8pt}
\noindent \textbf{(a)} No vertex $x \in Y$ is adjacent to two consecutive vertices on the cycle.

\vspace{8pt}
\noindent \textbf{(b)} If $x \in Y$ is adjacent to $u_i$ and $u_j$ then $u_{i+1}u_{j+1} \notin E(G)$.

\vspace{8pt}
\noindent \textbf{(c)} If $x \in Y$ is adjacent to $u_i$ and $u_j$ then no vertex $x' \in Y$ is adjacent to both
$u_{i+1}$ and $u_{j+2}$.

\vspace{8pt}
\noindent \textbf{(d)} Suppose $\alpha(G) = m - 1$ where $m \leq  \frac {n + 2}{2}$ and $\{x_1, x_2, ... ,x_{m -1} \} \subseteq Y$ is an $(m - 1)$-element independent set. Then, no member of this set is adjacent to $m -2$ or more vertices on the cycle (We have taken the liberty of making a slight correction to the inequality $m \leq \frac {n + 2}{2}$ of the original \cite{BoJaSh}, Lemma 5(d)).
\end{lemma}

\vspace{12pt}

\noindent The main results of this paper hinges on Lemma \ref{l4} that we prove next. 

\vspace{4pt}
\begin{lemma}
\label{l4}
\noindent A $C_n$ -free graph, where  $m \geq 7$ and $n \geq (m-3)(m-1)$ of order $(m-1)(n-1)$ with no independent set of order $m$  contains an isomorphic copy of  $(m-1)K_{n-1}$.
\end{lemma}

\noindent {\bf Proof.} We will prove this result using the principle of Mathematical Induction. By Lemma \ref{l2}, the result is true for $m=6$(\cite{JaNaSe}). In each of the two cases $n \geq (m-3)(m-1)+2$ and $n = (m-3)(m-1)+1$,  we consider  $G$ as a  graph on $(m-1)(n-1)$ vertices  satisfying $C_n \not \subseteq G$ and $\alpha(G) \leq m-1$. Since  $r(C_{n-1}, K_m)=(m-1)(n-2)+1 \leq (m-1)(n-1)$ (see \cite{BoJaSh,Ra}), there exists a cycle $C=(u_1, u_2, ... , u_{n-1}, u_1)$ of length $n - 1$ in $G$. In consistent with the notation of \cite{BoJaSh}, define $H$  as the induced subgraph of $G$ not containing the vertices of the cycle $C$. Then, $|V(C)|=n-1$ and $|V(H)|=(m-2)(n-1)$. 

\vspace{6pt}

\noindent Suppose  there exists an independent set $Y=\{y_1, y_2,..., y_{m-1} \}$ of size $m-1$ in $H$, so that $\alpha(G)=m-1$. From  Lemma \ref{l3}(d) (as $m\leq \frac{n+2}{2}$), it follows that no vertex of $Y$ is adjacent to $m-2$ or more vertices of $C$. Thus, $\left| E(Y,V(C)) \right| \le (m-1)(m-3)$.  For ease of reference, we define such a graph structure as the \textbf{Standard Configuration ($n$)}.

\vspace{10pt}

\noindent \textbf{Case 1: $n \geq  (m-3)(m-1)+2$}

\vspace{6pt}
\noindent Now, $\left| E(Y,V(C)) \right| \le (m-3(m-1)$ and $V(C)=(m-3)(m-1)+1$. Thus, there exists a vertex $x\in V(C)$ adjacent to no vertex of $Y$. This gives, an independent set $Y \cup \{x\}$ of size $m$, a contradiction.

\vspace{50pt}
\noindent \textbf{Case 2: $n = (m-3)(m-1)+1$}

\vspace{6pt}
\noindent By Lemma \ref{l3}, it follows that  $\left| E(Y,V(C)) \right| \le (m-1)(m-3)=n-1$. However if $\left| E(Y,V(C)) \right| \le (m-1)(m-3)$, we get a contradiction as in case 1. Thus, we get $\left| E(Y,V(C)) \right| = (m-1)(m-3)$.  Hence,  for  $1\leq i\leq m-1$, $\left|\Gamma (y_i) \cap V(C)\right| =m-3$ and for each $1\leq j  < j' \leq m-1$, $\Gamma (y_j) \cap \Gamma (y_{j'}) \cap V(C)=\phi$.

\vspace{10pt}

\noindent By Lemma \ref{l1}, as $\delta(G) \geq n-2$, we get that $\left|\Gamma (y_1) \cap V(H \setminus Y)\right| \geq n-2-(m-3)=(m-3)(m-2)-1$. Since $r(P_{m-3},K_m)=(m-4)(m-1)+1$ and $\alpha(G) < m$, we get that each of $G[\Gamma (y_1) \cap V(H \setminus Y)]$ and $G[\Gamma (y_2) \cap V(H \setminus Y)]$ contains a $P_{m-3}$. Thus, $P_{m-3} \subseteq \Gamma (y_1) \cap V(H \setminus Y)$, where $P_{m-3}$ is induced by  $\{x_1,x_2,...,x_{m-3}\}$ such that $(x_i,x_{i+1}) \in E(G)$ ($1 \leq i \leq {m-2}$) and $P_{m-3} \subseteq \Gamma (y_2) \cap V(H \setminus Y)$,  where $P_{m-3}$ is induced by  $\{x'_1,x'_2,...,x'_{m-3}\}$ such that $(x'_i,x'_{i+1}) \in E(G)$ ($1 \leq i \leq {m-2}$).

\vspace{10pt}

\noindent Suppose that $x_1$ is not adjacent to any vertex of $\{y_2,...,y_{m-1}\}$ and $x'_1$ is not adjacent to any vertex of $\{y_1,y_3,...,y_{m-1}\}$. Re-order the vertices of the cycle such that $y_1\in Y$ is adjacent to $u_1$. In this ordering, let $y_1$ be also adjacent to $u_t$ where $2 \leq t \leq (m-3)(m-1)$.

\begin{center}

\definecolor{uuuuuu}{rgb}{0.26666666666666666,0.26666666666666666,0.26666666666666666}
\definecolor{ffffff}{rgb}{1.0,1.0,1.0}
\begin{tikzpicture}[line cap=round,line join=round,>=triangle 45,x=1.089126559714795cm,y=1.2880658436213994cm]
\clip(-5.699999999999999,5.219999999999996) rectangle (5.520000000000013,10.07999999999999);
\fill[color=ffffff,fill=ffffff,fill opacity=0.1] (-3.02,5.62) -- (-2.44,5.58) -- (-1.8938741688442602,5.77936543467826) -- (-1.4760526241370269,6.183624209200794) -- (-1.2587805067477513,6.722876308876984) -- (-1.279626139628259,7.303880146721288) -- (-1.5349851234768812,7.826174880916263) -- (-1.980703570122651,8.199451007020894) -- (-2.5397127109155475,8.359165691724547) -- (-3.1153547869738647,8.2777028150796) -- (-3.608096053803033,7.969148048524964) -- (-3.9327370699544106,7.4868533143299905) -- (-4.033144454464381,6.9142117535470975) -- (-3.891956858378718,6.350238294320912) -- (-3.533586899709677,5.892449071289668) -- cycle;
\fill[color=ffffff,fill=ffffff,fill opacity=0.1] (1.08,6.36) -- (2.14,6.34) -- (2.486579144363348,7.341939567385363) -- (1.6407768353717538,7.9811722747028835) -- (0.7714631162884599,7.374300247160362) -- cycle;
\draw (-4.8199999999999985,1.820000000000001) node[anchor=north west] {$v_{1,3}$};
\draw (-0.6399999999999939,1.6400000000000012) node[anchor=north west] {$v_{ 2,3}$};
\draw (6.8212102632969694E-15,0.24000000000000307) node[anchor=north west] {$v_{2,1}$};
\draw (0.42000000000000726,7.319999999999993) node[anchor=north west] {$y_{1}$};
\draw (1.020000000000008,8.379999999999992) node[anchor=north west] {$y_{2}$};
\draw (1.6200000000000085,8.479999999999992) node[anchor=north west] {$y_{3}$};
\draw (2.180000000000009,8.259999999999993) node[anchor=north west] {$y_{4}$};
\draw (-1.3599999999999945,6.359999999999995) node[anchor=north west] {$u_{1}$};
\draw (6.8212102632969694E-15,0.24000000000000307) node[anchor=north west] {$u_{4}$};
\draw (-1.7799999999999951,8.619999999999992) node[anchor=north west] {$u_t$};
\draw [color=ffffff] (-3.02,5.62)-- (-2.44,5.58);
\draw [color=ffffff] (-2.44,5.58)-- (-1.8938741688442602,5.77936543467826);
\draw [color=ffffff] (-1.8938741688442602,5.77936543467826)-- (-1.4760526241370269,6.183624209200794);
\draw [color=ffffff] (-1.4760526241370269,6.183624209200794)-- (-1.2587805067477513,6.722876308876984);
\draw [color=ffffff] (-1.2587805067477513,6.722876308876984)-- (-1.279626139628259,7.303880146721288);
\draw [color=ffffff] (-1.279626139628259,7.303880146721288)-- (-1.5349851234768812,7.826174880916263);
\draw [color=ffffff] (-1.5349851234768812,7.826174880916263)-- (-1.980703570122651,8.199451007020894);
\draw [color=ffffff] (-1.980703570122651,8.199451007020894)-- (-2.5397127109155475,8.359165691724547);
\draw [color=ffffff] (-2.5397127109155475,8.359165691724547)-- (-3.1153547869738647,8.2777028150796);
\draw [color=ffffff] (-3.1153547869738647,8.2777028150796)-- (-3.608096053803033,7.969148048524964);
\draw [color=ffffff] (-3.608096053803033,7.969148048524964)-- (-3.9327370699544106,7.4868533143299905);
\draw [color=ffffff] (-3.9327370699544106,7.4868533143299905)-- (-4.033144454464381,6.9142117535470975);
\draw [color=ffffff] (-4.033144454464381,6.9142117535470975)-- (-3.891956858378718,6.350238294320912);
\draw [color=ffffff] (-3.891956858378718,6.350238294320912)-- (-3.533586899709677,5.892449071289668);
\draw [color=ffffff] (-3.533586899709677,5.892449071289668)-- (-3.02,5.62);
\draw [color=ffffff] (1.08,6.36)-- (2.14,6.34);
\draw [color=ffffff] (2.14,6.34)-- (2.486579144363348,7.341939567385363);
\draw [color=ffffff] (2.486579144363348,7.341939567385363)-- (1.6407768353717538,7.9811722747028835);
\draw [color=ffffff] (1.6407768353717538,7.9811722747028835)-- (0.7714631162884599,7.374300247160362);
\draw [color=ffffff] (0.7714631162884599,7.374300247160362)-- (1.08,6.36);
\draw (6.8212102632969694E-15,0.24000000000000307) node[anchor=north west] {$y_{2}$};
\draw (2.5600000000000094,7.779999999999993) node[anchor=north west] {$y_{5}$};
\draw (-1.4760526241370269,6.183624209200794)-- (-1.2587805067477513,6.722876308876984);
\draw (-1.2587805067477513,6.722876308876984)-- (-1.279626139628259,7.303880146721288);
\draw (-1.279626139628259,7.303880146721288)-- (-1.5349851234768812,7.826174880916263);
\draw (-1.5349851234768812,7.826174880916263)-- (-1.980703570122651,8.199451007020894);
\draw (-1.980703570122651,8.199451007020894)-- (-2.5397127109155475,8.359165691724547);
\draw (-2.5397127109155475,8.359165691724547)-- (-3.1153547869738647,8.2777028150796);
\draw (-3.1153547869738647,8.2777028150796)-- (-3.608096053803033,7.969148048524964);
\draw (-3.608096053803033,7.969148048524964)-- (-3.9327370699544106,7.4868533143299905);
\draw (-3.9327370699544106,7.4868533143299905)-- (-4.033144454464381,6.9142117535470975);
\draw (-4.033144454464381,6.9142117535470975)-- (-3.891956858378718,6.350238294320912);
\draw (-3.891956858378718,6.350238294320912)-- (-3.533586899709677,5.892449071289668);
\draw (-3.533586899709677,5.892449071289668)-- (-3.02,5.62);
\draw (-3.02,5.62)-- (-2.44,5.58);
\draw (-2.44,5.58)-- (-1.8938741688442602,5.77936543467826);
\draw (-1.8938741688442602,5.77936543467826)-- (-1.4760526241370269,6.183624209200794);
\draw (0.7714631162884599,7.374300247160362)-- (-1.4760526241370269,6.183624209200794);
\draw [->] (0.7714631162884599,7.374300247160362) -- (-0.9800000000000001,7.96);
\draw (-0.78,9.26)-- (-0.1,9.28);
\draw (-0.1,9.28)-- (1.4,9.28);
\draw (-0.25999999999999346,9.77999999999999) node[anchor=north west] {$x_2$};
\draw (-0.9799999999999942,9.77999999999999) node[anchor=north west] {$x_1$};
\draw (1.3800000000000083,9.759999999999991) node[anchor=north west] {$x_{m-3}$};
\draw (0.7714631162884599,7.374300247160362)-- (-0.78,9.26);
\draw (0.7714631162884599,7.374300247160362)-- (-0.1,9.28);
\draw (0.7714631162884599,7.374300247160362)-- (1.4,9.28);
\draw (-2.9999999999999964,7.199999999999994) node[anchor=north west] {$C$};
\draw (1.3600000000000083,7.399999999999993) node[anchor=north west] {$Y$};
\draw [->] (-0.78,9.26) -- (-1.24,8.14);
\begin{scriptsize}
\draw [fill=black] (-3.02,5.62) circle (1.5pt);
\draw [fill=black] (-2.44,5.58) circle (1.5pt);
\draw [fill=uuuuuu] (-1.8938741688442602,5.77936543467826) circle (1.5pt);
\draw [fill=uuuuuu] (-1.4760526241370269,6.183624209200794) circle (1.5pt);
\draw [fill=uuuuuu] (-1.2587805067477513,6.722876308876984) circle (1.5pt);
\draw [fill=uuuuuu] (-1.279626139628259,7.303880146721288) circle (1.5pt);
\draw [fill=uuuuuu] (-1.5349851234768812,7.826174880916263) circle (1.5pt);
\draw [fill=uuuuuu] (-1.980703570122651,8.199451007020894) circle (1.5pt);
\draw [fill=uuuuuu] (-2.5397127109155475,8.359165691724547) circle (1.5pt);
\draw [fill=uuuuuu] (-3.1153547869738647,8.2777028150796) circle (1.5pt);
\draw [fill=uuuuuu] (-3.608096053803033,7.969148048524964) circle (1.5pt);
\draw [fill=uuuuuu] (-3.9327370699544106,7.4868533143299905) circle (1.5pt);
\draw [fill=uuuuuu] (-4.033144454464381,6.9142117535470975) circle (1.5pt);
\draw [fill=uuuuuu] (-3.891956858378718,6.350238294320912) circle (1.5pt);
\draw [fill=uuuuuu] (-3.533586899709677,5.892449071289668) circle (1.5pt);
\draw [fill=black] (2.14,6.34) circle (1.5pt);
\draw [fill=uuuuuu] (2.486579144363348,7.341939567385363) circle (1.5pt);
\draw [fill=uuuuuu] (1.6407768353717538,7.9811722747028835) circle (1.5pt);
\draw [fill=black] (0.7714631162884599,7.374300247160362) circle (1.5pt);
\draw [fill=black] (-0.78,9.26) circle (1.5pt);
\draw [fill=black] (1.4,9.28) circle (1.5pt);
\draw [fill=black] (-0.1,9.28) circle (1.5pt);
\draw [fill=black] (1.08,7.8) circle (1.5pt);
\draw [fill=black] (2.12,7.8) circle (1.5pt);
\draw [fill=black] (2.42,6.78) circle (1.5pt);
\draw [fill=black] (0.46,9.58) circle (1.0pt);
\draw [fill=black] (0.82,9.58) circle (1.0pt);
\draw [fill=black] (1.16,9.56) circle (1.0pt);
\end{scriptsize}
\end{tikzpicture}
\end{center}
\begin{center}
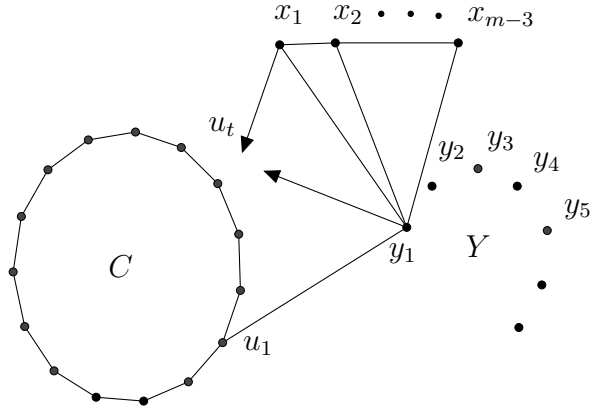
Figure 1. Configuration for $n-1=(m-3)(m-1)$
\end{center}

\noindent By Lemma \ref{l2}(a), $t \neq 2$. In order to avoid an independent set of size $m$, induced by $\{x,u_t,y_2,y_3,...,y_{m-1}\}$, we get that $(x_1,u_t) \in E(G)$. However, $t \neq i$ where $3 \leq i \leq {m-1}$, since otherwise  $t=i$ and then we get a  $C_{n}$ comprising \[(u_1,y_1,x_{i-2},x_{i-3},...,x_{1}, u_t,...,u_{(m-3)(m-1)},u_1).\] 

\vspace{8pt}
\noindent Thus, any pair of vertices adjacent to $y_1$ in $C$ cannot be separated by a path of length 1, 2,..., $(m-2)$  along $C$.  Hence, $\Gamma (y_1) \cap C =\{u_1,u_m,u_{2m-1},....,u_{(m-4)(m-1)+1}\}$, as $\frac{n-1}{m-3}=({m-1})$. In this scenario, we use the prerogative that $(y_2,u_2)\in E(G)$. Then,  by the same argument $\Gamma (y_2) \cap C =\{u_2,u_{m+1},u_{2m},....,u_{(m-4)(m-1)+2} \}$. But by Lemma \ref{l3}(b),  $(u_2,u_{m+1})\notin E(G)$. Henceforth, we will get that $\{u_2,u_{m+1},y_1,y_3,y_4,...,y_{m-1}\}$  is an independent set of size $m$, a contradiction. 
\vspace{10pt}

\noindent This implies that  $x_1$ is adjacent to some vertex of $\{y_2,y_3,...,y_{m-1}\}$ or  $x'_1$ is adjacent to some vertex of $\{y_1,y_3,y_4,...,y_{m-1}\}$. Therefore, without loss of generality, by relabeling if necessary, we may assume that $y_1$ is adjacent to $\{x_1,x_2,x_3,...,x_{m-3} \}\subseteq V(H \setminus Y)$ where $(x_i,x_{i+1}) \in E(G)$, $1 \leq i \leq {m-4}$. Let $y_1$ be also adjacent to $X'=\{u_{i_1},u_{i_2},u_{i_3},...,u_{i_{m-3}}\} \subseteq V(C)$ where $1=i_1<i_2<i_3<...<i_{m-3} \leq n-3$, $X'$ induces a  $K_{m-3}$  and $y_2$ is adjacent to $x_1$. Define $S_{q} = \{u_i : i_q \leq i \leq i_{(q+1)} \}$, whenever $1 \leq q \leq {m-4}$ and $S_{m-3} = \{u_{i_{m-3}},....,u_{n-1},u_1\}$.

\vspace{14pt}
\noindent \textbf{Remark:} As a result of Lemma \ref{l3}, $\{u_{i_1+1},u_{i_2+1},u_{i_3+1}, u_{i_4+1},...,u_{i_{(m-4)}+1},u_{i_{(m-3)}+1} \}$ and $\{u_{i_1+1},u_{i_2+1},u_{i_3+1}, u_{i_4+1},...,u_{i_{(m-4)}+1}, u_{i_{(m-3)}+2} \}$  are independent sets of size $m-3$. Moreover, $y_2$  is adjacent to exactly one vertex of \[\{u_{i_1+1},u_{i_2+1},u_{i_3+1},u_{i_4+1},...,u_{i_{(m-4)}+1},u_{i_{(m-3)}+1}, u_{i_{(m-3)}+2}\}\] and no other vertex of $\{y_1, y_3, y_4,...,y_{m-1}\}$ is adjacent to that vertex.

\vspace{14pt}

\noindent \textbf{Claim 1:}
\label{l5}

\vspace{4pt}
\noindent \textbf{a)} There exists some $1 \leq q \leq {m-3}$ such that $\left|S_q \right|\geq m$. 

\noindent \textbf{b)} There cannot exist $1 \leq q_1 < q_2 < q_3 \leq {m-3}$ such that $\left|S_{q_1} \right|=\left|S_{q_2} \right|=\left|S_{q_3} \right|=3$. 

\noindent \textbf{c)} There cannot exist  $1 \leq q_1, q_2 \leq {m-3}$ such that $\left|S_{q_1} \right|=3$ and $\left|S_{q_2} \right|=4$. 

\noindent \textbf{d)} There cannot exist  $1 \leq q_1\leq {m-3}$ such that $\left|S_{q_1} \right|=5$. 

\noindent \textbf{e)} If for any $1 \leq q\leq {m-3}$,  $\left|S_{q} \right|\geq 5$ then $\left|S_{q} \right|\geq m$.

\vspace{10pt}
\noindent {\bf Proof of Claim 1.} \textbf{a)} Suppose that (a) is false. That is, for all $1 \leq q \leq {m-3}$,  $\left|S_q \right| \leq m-1$. Notice that, $\bigcup_{i=1}^{m-3} S_i = V(C)$. Then \[ \left|\bigcup_{i=1}^{m-3} S_i\right| =\left|V(C)\right|=(m-3)(m-1)\].

\noindent  Note that any pair of  $S_i$'s are disjoint unless they are consecutive $S_i$'s and in such a case there is exactly one element in common. Also any intersection of three or more $S_i$'s will have empty intersection.  Under these conditions,  \[ \left|\bigcup_{i=1}^{m-3} S_i\right| =\sum_{i=1}^{m-3} \left|S_i\right|-\sum_{i=1,j=1}^{m-3} \left|S_i \cap S_i \right| \leq  (m-3)(m-1) -(m-3) =  (m-3)(m-2)  \]. 

\noindent This gives us $m \leq  2$, a contradiction. Hence, (a) is true.

\vspace{8pt}
\noindent \textbf{b)} Suppose that (b) is false. Since the position of the three $S_i$'s do not play a role, without loss of generality, assume that  $\left|S_{1} \right|=\left|S_{2} \right|=\left|S_{3} \right|=3$. That is, $y_1$ is adjacent to $\{u_1, u_3, u_5, u_7, u_{i_5},...,u_{i_{(m-3)}} \}$. By the above remark, $y_2$ is adjacent to exactly one vertex of $\{u_2, u_4, u_6, u_{i_4+1},...,u_{i_{(m-3)}+1}, u_{i_{(m-3)}+2} \}$. If $y_2$ is adjacent to $u_2$, we get a $C_n$ consisting of  $(u_1,y_1,x_{1},y_{2},u_2, u_3,u_5,u_7,u_8,...,u_{(m-3)(m-1)},u_1)$, a contradiction. Similar argument follows when $y_2$ is adjacent to any vertex of $\{u_4,u_6,u_{i_4+1},...,u_{i_{(m-3)}+1} \}$.  Finally, for the remaining possibility that $y_2$ is adjacent to $u_{i_{(m-3)}+2}$, gives a $C_n$ consisting of  $(u_1,...,u_5,u_7,u_8,u_9,....,u_{i_{(m-3)}},y_1,x_1,y_2,u_{i_{(m-3)}+2},...,u_{(m-3)(m-1)},u_1)$, a contradiction. 

\vspace{8pt}
\noindent \textbf{c)} and \textbf{d)} The argument used to prove $(b)$ was purely based on the total number of interior points contained in the three $S_i$'s. Since there are a total of three interior points, one from first $S_i$ and two from one $S_i$ mentioned in part $(c)$  and there are a total of three interior points in the $S_i$ mentioned  in part $(d)$,  the results of $(c)$ and $(d)$ follow through by the same argument used in $(b)$.

\vspace{8pt}
\noindent \textbf{e)} Suppose that for all $1 \leq q\leq {m-3}$,  $\left|S_{q} \right|\geq 5$. 
and that there exists $1 \leq q_1\leq {m-3}$ such that $5\leq \left|S_{q_1} \right| \leq m-1$. By the above remark, $y_2$ is adjacent to exactly one vertex of $\{u_{i_1+1}, u_{i_2+1},...,u_{i_{(m-3)}+1}, u_{i_{(m-3)}+2} \}$. If $y_2$ is adjacent to $u_{i_j+1}$ where $i_j<q_1$, we get a $C_n$ consisting of  $(u_1,...,u_{i_j},y_1,x_{m-4},...,x_1,y_2,u_{i_j+1},...,u_{q_1},u_{q_1+1},...,u_{n-1},u_1)$, a contradiction. Next if $y_2$ is adjacent to $u_{i_j+1}$ where $q_1<i_j$, we get a $C_n$ consisting of  $(u_1,...,u_{q_1},u_{q_1+1},...,u_{i_j},y_1,x_{m-4},...,x_1,y_2,u_{i_j+1},...,u_{n-1},u_1)$, a contradiction. Also if $y_2$ is adjacent to $u_{i_j+2}$ where $q_1<i_j$, we get a $C_n$ consisting of  $(u_1,...,u_{q_1},u_{q_1+1},...,u_{i_j},y_1,x_{m-3},...,x_1,y_2,u_{i_j+1},...,u_{n-1},u_1)$, a contradiction. All the remaining possibilities will be similar to one of the three possibilities we have already considered and therefore the result will follow likewise.

\vspace{8pt}
\noindent Next we continue with the proof of case 2 of  Lemma \ref{l4}. By the claim we get that:  $\left|S_{q} \right|\geq m$ whenever $1 \leq q\leq {m-3}$ except if there exists,  $1 \leq q_1 \leq {m-3}$ such that $\left|S_{q_1} \right|=3$ or else there exists, $1 \leq q_1<q_2 \leq {m-3}$ such that  $\left|S_{q_1} \right|=\left|S_{q_2} \right|=3$ or else there exists, $1 \leq q_1 \leq {m-3}$ such that  $\left|S_{q_1} \right|=4$. However in all these three possibilities all the other $S_{q}$'s will satisfy $\left|S_{q} \right|\geq m$. \noindent In the first possibility, there exists,  $1 \leq q_1 \leq {m-3}$ such that $\left|S_{q_1} \right|=3$ and all the other $S_{q}$'s will satisfy $\left|S_{q} \right|\geq m$ we get that  \[ (m-3)+ 1+(m-2)(m-4) \geq \left|\bigcup_{i=1}^{m-3} S_i\right| = \left|V(C)\right|=(m-3)(m-1)\]

\noindent This gives us $m \leq  3$, a contradiction. 

\vspace{8pt}
\noindent In the second possibility, there exists,  $1 \leq q_1<q_2 \leq {m-3}$ such that  $\left|S_{q_1} \right|=\left|S_{q_2} \right|=3$ and all the other $S_{q}$'s will satisfy $\left|S_{q} \right|\geq m$ we get that  \[ (m-3)+ 2+(m-2)(m-5) \geq \left|\bigcup_{i=1}^{m-3} S_i\right| = \left|V(C)\right|=(m-3)(m-1)\]

\noindent This gives us $m \leq  3$, a contradiction. 

\vspace{8pt}
\noindent In the third possibility, there exists,  $1 \leq q_1 \leq {m-3}$ such that  $\left|S_{q_1} \right|=4$ and all the other $S_{q}$'s will satisfy $\left|S_{q} \right|\geq m$ we get that  \[ (m-2)+(m-2)(m-4) \geq \left|\bigcup_{i=1}^{m-3} S_i\right| = \left|V(C)\right|=(m-3)(m-1)\]

\noindent This gives us $m \leq  3$, a contradiction.  Thus we can conclude that $\left|S_{q} \right|\geq m$ whenever $1 \leq q\leq {m-3}$. 

\vspace{8pt}
\noindent Thus, any pair of vertices adjacent to $y_1$ (or $y_2$) in $C$ cannot be separated by a path of length 1, 2,..., $(m-2)$  along $C$.  We find $\Gamma (y_1) \cap C =\{u_1,u_m,u_{2m-1},....,u_{1+(m-4)(m-1)}\}$, as $\frac{n-1}{m-3}=({m-1})$. In this scenario, if we consider the possibility that $y_2$ is adjacent to $u_2$, we get a $C_n$ consisting of  $(u_1,u_2,y_1,x_1,...,x_{m-4},y_2,u_m,...,u_{(m-3)(m-1)},u_1)$, a contradiction. Also if we consider the possibility that $y_2$ is adjacent to $u_{(m-4)(m-1)+3}$, we get a $C_n$ consisting of  $(u_1,u_2,...,u_{(m-4)(m-1)+3},y_1,x_1,...,x_{m-5},y_2,,u_{(m-3)(m-1)},u_1)$, a contradiction. All the remaining $m-5$ possibilities will be similar to one of the two possibilities we have already considered and therefore the result will follow likewise.
 
\vspace{12pt} 
\noindent  Since this is impossible this concludes the proof of case 2 of Lemma \ref{l4}.

\vspace{12pt}
\noindent Having proved that  $H$ cannot have an independent set of size $m-1$ in both cases $n = (m-3)(m-1)+2$ and $n = (m-3)(m-1)+1$,  we next continue with the proof of Lemma \ref{l4}.

\vspace{12pt}
\noindent Since, $H$  satisfies all conditions of the induction hypothesis, $H$ contains an isomorphic copy of $(m-2)K_{n-1}$.

\vspace{12pt}
\noindent Next we show that $V(C_{n-1})$ induced a $K_{n-1}$. Suppose that there exists two vertices  of $V(C)$, say $v$ and $w$, such that $(v,w) \not \in E(G)$. In order to avoid a $C_n$ both $v$ and $w$ will have to be adjacent to at most one vertex of each of the $m-2$ copies of $K_{n-1}$ in $H$. Moreover, any vertex of any copy of $K_{n-1}$ in $H$ will have to be adjacent to at most one vertex of another copy of a $K_{n-1}$ in $H$. Thus, each copy of a $K_{n-1}$  will have at most $m-2$ vertices adjacent to some vertex outside that of $K_{n-1}$, in $V(H) \cup \{v,w\}$. Since $(n-1)-(m-2) \geq 1$, we can select $x_1$ in the first $K_{n-1}$, $x_2$ in the second $K_{n-1}$, $ $... $ $ and likewise $x_{m-2}$ in the $(m-2)^{\textsl{th}}$ $K_{n-1}$ such that $\{x_1, x_2,..., x_{m-2}\}$ is an independent set of size $m-2$ and no vertex of $\{x_1, x_2,..., x_{m-2}\}$ is adjacent to any vertex of $\{v, w\}$. Hence $\{x_1, x_2, ,...,x_{m-2}, v, w\}$ is an independent set of size $m$, a contradiction. Therefore, we get that any two pairs of vertices of $V(C)$ are  connected by an edge. Hence, $G[V(C_{n-1})]=K_{n-1}$ as required. This $K_{n-1}$ along with the $(m-2)K_{n-1}$ contained in $H$  gives the required $(m-1)K_{n-1}$.

\section{MAIN RESULT}

\begin{theorem}
\label{t1}

\noindent If $m \geq 7$ and $n \geq (m-3)(m-1)$, then $r_*(C_n, K_m) =(m-2)(n-1)+2$. 
\end{theorem}

\noindent {\bf Proof.}\vspace{6pt}
\noindent  To find a lower bound for $r_*(C_n, K_m)$, color the graph $K_{(m-1)(n-1)+1}\setminus K_{1,n-2}$, such that the red graph consists of a $(m-2)K_{n-1} \cup ( K_{n-1} \sqcup K_{1,1} )$ as illustrated in the following figure.  

\begin{center}

\begin{tikzpicture}[line cap=round,line join=round,>=triangle 45,x=1.0cm,y=1.0cm]
\clip(-5.994545454545459,-2.069090909090903) rectangle (9.460000000000006,7.330909090909107);
\draw (1.7872727272727287,6.640000000000016) node[anchor=north west] {red degree 1};
\draw (1.7690909090909106,7.203636363636379) node[anchor=north west] {blue degree};
\draw (3.9690909090909123,7.258181818181834) node[anchor=north west] {$(m-2)(n-1)$};
\draw (6.569090909090914,-1.2509090909090843) node[anchor=north west] {$K_{n-1}$};
\draw (5.823636363636368,-1.2509090909090843) node[anchor=north west] {red};
\draw (-0.5581818181818183,-1.2509090909090843) node[anchor=north west] {$K_{n-1}$};
\draw (-1.3036363636363644,-1.2509090909090843) node[anchor=north west] {red};
\draw (2.823636363636366,4.312727272727286) node[anchor=north west] {$K_{n-1}$};
\draw (1.9872727272727289,4.312727272727286) node[anchor=north west] {red};
\draw (2.914545454545458,5.403636363636376)-- (2.92,3.3);
\draw(2.92,3.3) circle (0.4cm);
\draw(0.2781818181818182,2.367272727272727) circle (0.4cm);
\draw(5.514545454545458,2.24) circle (0.4cm);
\draw(-0.14,1.24) circle (0.4cm);
\draw(-0.01272727272727264,0.0036363636363640414) circle (0.4cm);
\draw(0.8054545454545451,-1.16) circle (0.4cm);
\draw(4.878181818181818,-1.250909090909091) circle (0.4cm);
\draw(5.678181818181818,-0.01454545454545416) circle (0.4cm);
\draw(5.8963636363636365,1.2218181818181821) circle (0.4cm);
\draw(1.4054545454545453,3.1854545454545455) circle (0.4cm);
\draw(4.423636363636367,2.9854545454545454) circle (0.4cm);
\draw [line width=2.0pt,dash pattern=on 4pt off 4pt] (2.4054545454545475,3.3490909090909207)-- (1.9690909090909108,3.276363636363648);
\draw [line width=2.0pt,dash pattern=on 4pt off 4pt] (0.8963636363636373,2.9672727272727384)-- (0.5327272727272733,2.730909090909102);
\draw [line width=2.0pt,dash pattern=on 4pt off 4pt] (-0.1036363636363635,0.749090909090918)-- (-0.060385360974644114,0.400787099031133);
\draw [line width=2.0pt,dash pattern=on 4pt off 4pt] (0.172299129950738,-0.3509975748298164)-- (0.49636363636363695,-0.76);
\draw [line width=2.0pt,dash pattern=on 4pt off 4pt] (5.478181818181822,-0.5054545454545379)-- (5.26,-0.869090909090902);
\draw [line width=2.0pt,dash pattern=on 4pt off 4pt] (5.823636363636368,0.7127272727272815)-- (5.775196068196354,0.37351154551267796);
\draw [line width=2.0pt,dash pattern=on 4pt off 4pt] (5.660305193977666,1.8675028881177)-- (5.803715425441227,1.6109406676923304);
\draw [line width=2.0pt,dash pattern=on 4pt off 4pt] (4.9145454545454585,2.767272727272738)-- (5.264667435523693,2.552347523777216);
\draw [line width=2.0pt,dash pattern=on 4pt off 4pt] (3.319055951448792,3.3274653893709796)-- (4.069002425170189,3.170480948132561);
\draw [line width=2.0pt,dash pattern=on 4pt off 4pt] (0.7690909090909098,1.9854545454545558)-- (5.150909090909095,0.1127272727272809);
\draw [line width=2.0pt,dash pattern=on 4pt off 4pt] (0.7327272727272734,-0.0145454545454465)-- (5.150909090909095,0.0945454545454627);
\draw [line width=2.0pt,dash pattern=on 4pt off 4pt] (1.16909090909091,-0.6327272727272651)-- (3.9872727272727304,2.64);
\draw [line width=2.0pt,dash pattern=on 4pt off 4pt] (1.16909090909091,-0.6327272727272651)-- (2.8963636363636387,2.84);
\draw [line width=2.0pt,dash pattern=on 4pt off 4pt] (2.8963636363636387,2.84)-- (4.823636363636368,-0.6509090909090836);
\draw [line width=2.0pt,dash pattern=on 4pt off 4pt] (0.7327272727272734,-0.0145454545454465)-- (4.823636363636368,-0.6509090909090836);
\draw [line width=2.0pt,dash pattern=on 4pt off 4pt] (5.150909090909095,0.0945454545454627)-- (0.49636363636363695,1.2218181818181912);
\draw [line width=2.0pt,dash pattern=on 4pt off 4pt] (4.9145454545454585,2.1854545454545558)-- (0.7327272727272734,-0.0145454545454465);
\draw [line width=2.0pt,dash pattern=on 4pt off 4pt] (0.7327272727272734,-0.0145454545454465)-- (1.714545454545456,2.767272727272738);
\draw [line width=2.0pt,dash pattern=on 4pt off 4pt] (0.7327272727272734,-0.0145454545454465)-- (0.7690909090909098,1.9854545454545558);
\draw [line width=2.0pt,dash pattern=on 4pt off 4pt] (0.49636363636363695,1.2218181818181912)-- (1.714545454545456,2.767272727272738);
\draw [line width=2.0pt,dash pattern=on 4pt off 4pt] (1.714545454545456,2.767272727272738)-- (1.16909090909091,-0.6327272727272651);
\draw [line width=2.0pt,dash pattern=on 4pt off 4pt] (0.7690909090909098,1.9854545454545558)-- (4.9145454545454585,2.1854545454545558);
\draw [line width=2.0pt,dash pattern=on 4pt off 4pt] (5.150909090909095,0.0945454545454627)-- (1.714545454545456,2.767272727272738);
\draw [line width=2.0pt,dash pattern=on 4pt off 4pt] (4.823636363636368,-0.6509090909090836)-- (4.9145454545454585,2.1854545454545558);
\draw [line width=2.0pt,dash pattern=on 4pt off 4pt] (4.823636363636368,-0.6509090909090836)-- (4.005454545454545,2.712727272727273);
\draw [line width=2.0pt,dash pattern=on 4pt off 4pt] (4.823636363636368,-0.6509090909090836)-- (1.714545454545456,2.767272727272738);
\draw [line width=2.0pt,dash pattern=on 4pt off 4pt] (4.823636363636368,-0.6509090909090836)-- (0.7690909090909098,1.9854545454545558);
\draw [line width=2.0pt,dash pattern=on 4pt off 4pt] (4.823636363636368,-0.6509090909090836)-- (0.49636363636363695,1.2218181818181912);
\draw [line width=2.0pt,dash pattern=on 4pt off 4pt] (4.823636363636368,-0.6509090909090836)-- (1.16909090909091,-0.6327272727272651);
\draw [line width=2.0pt,dash pattern=on 4pt off 4pt] (1.16909090909091,-0.6327272727272651)-- (5.150909090909095,0.0945454545454627);
\draw [line width=2.0pt,dash pattern=on 4pt off 4pt] (1.16909090909091,-0.6327272727272651)-- (4.9145454545454585,2.1854545454545558);
\draw [line width=2.0pt,dash pattern=on 4pt off 4pt] (2.914545454545458,5.403636363636376)-- (1.6054545454545468,3.6763636363636483);
\draw [line width=2.0pt,dash pattern=on 4pt off 4pt] (2.914545454545458,5.403636363636376)-- (4.223636363636367,3.5309090909091028);
\draw [line width=2.0pt,dash pattern=on 4pt off 4pt] (2.914545454545458,5.403636363636376)-- (0.33272727272727315,3.021818181818193);
\draw [line width=2.0pt,dash pattern=on 4pt off 4pt] (2.914545454545458,5.403636363636376)-- (5.441818181818186,2.8763636363636476);
\draw [shift={(2.354048779089766,1.685041843402511)},line width=2.0pt,dash pattern=on 4pt off 4pt]  plot[domain=0.00736204058457097:1.421194380074314,variable=\t]({1.0*3.7605985865431335*cos(\t r)+-0.0*3.7605985865431335*sin(\t r)},{0.0*3.7605985865431335*cos(\t r)+1.0*3.7605985865431335*sin(\t r)});
\draw [shift={(3.3705460019157956,1.665968161291805)},line width=2.0pt,dash pattern=on 4pt off 4pt]  plot[domain=1.6921977178825343:3.1291741749836173,variable=\t]({1.0*3.765381798705111*cos(\t r)+-0.0*3.765381798705111*sin(\t r)},{0.0*3.765381798705111*cos(\t r)+1.0*3.765381798705111*sin(\t r)});
\draw [shift={(2.528890479599143,2.0849906943450356)},line width=2.0pt,dash pattern=on 4pt off 4pt]  plot[domain=1.4551066964422603:3.6689506401506793,variable=\t]({1.0*3.340978724566652*cos(\t r)+-0.0*3.340978724566652*sin(\t r)},{0.0*3.340978724566652*cos(\t r)+1.0*3.340978724566652*sin(\t r)});
\draw [shift={(3.1486041479248255,2.023510418439554)},line width=2.0pt,dash pattern=on 4pt off 4pt]  plot[domain=-0.5045986730445833:1.639931521788343,variable=\t]({1.0*3.388220016076165*cos(\t r)+-0.0*3.388220016076165*sin(\t r)},{0.0*3.388220016076165*cos(\t r)+1.0*3.388220016076165*sin(\t r)});
\draw [shift={(2.0687033067973077,2.0450123865724077)},line width=2.0pt,dash pattern=on 4pt off 4pt]  plot[domain=1.3240848903214686:4.1129063015334975,variable=\t]({1.0*3.4634959157210203*cos(\t r)+-0.0*3.4634959157210203*sin(\t r)},{0.0*3.4634959157210203*cos(\t r)+1.0*3.4634959157210203*sin(\t r)});
\draw [shift={(3.5911771648736313,1.9642161956831503)},line width=2.0pt,dash pattern=on 4pt off 4pt]  plot[domain=-0.9961844730609712:1.7650441708367315,variable=\t]({1.0*3.5053447138826477*cos(\t r)+-0.0*3.5053447138826477*sin(\t r)},{0.0*3.5053447138826477*cos(\t r)+1.0*3.5053447138826477*sin(\t r)});
\draw [line width=2.0pt,dash pattern=on 4pt off 4pt] (4.9145454545454585,2.1854545454545558)-- (5.150909090909095,0.0945454545454627);
\draw [line width=2.0pt,dash pattern=on 4pt off 4pt] (4.005454545454545,2.712727272727273)-- (3.9872727272727304,2.64);
\draw [line width=2.0pt,dash pattern=on 4pt off 4pt] (5.223636363636367,1.276363636363646)-- (4.005454545454545,2.712727272727273);
\draw [line width=2.0pt,dash pattern=on 4pt off 4pt] (5.223636363636367,1.276363636363646)-- (4.823636363636368,-0.6509090909090836);
\draw [line width=2.0pt,dash pattern=on 4pt off 4pt] (5.223636363636367,1.276363636363646)-- (0.7327272727272734,-0.0145454545454465);
\draw [line width=2.0pt,dash pattern=on 4pt off 4pt] (5.223636363636367,1.276363636363646)-- (0.49636363636363695,1.2218181818181912);
\draw [line width=2.0pt,dash pattern=on 4pt off 4pt] (5.223636363636367,1.276363636363646)-- (0.7690909090909098,1.9854545454545558);
\draw [line width=2.0pt,dash pattern=on 4pt off 4pt] (5.223636363636367,1.276363636363646)-- (1.714545454545456,2.767272727272738);
\draw [line width=2.0pt,dash pattern=on 4pt off 4pt] (5.223636363636367,1.276363636363646)-- (2.8963636363636387,2.84);
\draw [line width=2.0pt,dash pattern=on 4pt off 4pt] (0.49636363636363695,1.2218181818181912)-- (2.8963636363636387,2.84);
\draw [line width=2.0pt,dash pattern=on 4pt off 4pt] (0.7327272727272734,-0.0145454545454465)-- (2.8963636363636387,2.84);
\draw [line width=2.0pt,dash pattern=on 4pt off 4pt] (0.7690909090909098,1.9854545454545558)-- (2.8963636363636387,2.84);
\draw [line width=2.0pt,dash pattern=on 4pt off 4pt] (2.8963636363636387,2.84)-- (4.9145454545454585,2.1854545454545558);
\draw [line width=2.0pt,dash pattern=on 4pt off 4pt] (1.714545454545456,2.767272727272738)-- (4.9145454545454585,2.1854545454545558);
\draw [line width=2.0pt,dash pattern=on 4pt off 4pt] (4.005454545454545,2.712727272727273)-- (0.7690909090909098,1.9854545454545558);
\draw [line width=2.0pt,dash pattern=on 4pt off 4pt] (4.005454545454545,2.712727272727273)-- (0.49636363636363695,1.2218181818181912);
\draw [line width=2.0pt,dash pattern=on 4pt off 4pt] (0.49636363636363695,1.2218181818181912)-- (4.9145454545454585,2.1854545454545558);
\draw [line width=2.0pt,dash pattern=on 4pt off 4pt] (0.7690909090909098,1.9854545454545558)-- (1.16909090909091,-0.6327272727272651);
\draw [line width=2.0pt,dash pattern=on 4pt off 4pt] (2.8963636363636387,2.84)-- (5.150909090909095,0.0945454545454627);
\draw [line width=2.0pt,dash pattern=on 4pt off 4pt] (4.005454545454545,2.712727272727273)-- (5.150909090909095,0.0945454545454627);
\draw [line width=2.0pt,dash pattern=on 4pt off 4pt] (5.223636363636367,1.276363636363646)-- (1.16909090909091,-0.6327272727272651);
\draw [line width=2.0pt,dash pattern=on 4pt off 4pt] (0.11266064046134114,2.003126136287687)-- (-0.10678180805850461,1.638618303297952);
\draw [line width=2.0pt,dash pattern=on 4pt off 4pt] (4.005454545454541,2.712727272727289)-- (1.714545454545456,2.767272727272738);
\draw(2.914545454545458,5.403636363636376) circle (0.1cm);
\begin{scriptsize}
\draw [fill=black] (1.918181818181818,-1.281818181818182) circle (2.0pt);
\draw [fill=black] (2.358181818181818,-1.281818181818182) circle (2.0pt);
\draw [fill=black] (2.778181818181818,-1.301818181818182) circle (2.0pt);
\draw [fill=black] (3.118181818181818,-1.301818181818182) circle (0.5pt);
\draw [fill=black] (3.378181818181818,-1.301818181818182) circle (0.5pt);
\draw [fill=black] (3.6963636363636367,-1.2872727272727273) circle (2.0pt);
\draw [fill=black] (2.914545454545458,5.403636363636376) circle (2.0pt);
\end{scriptsize}
\end{tikzpicture}
\end{center}
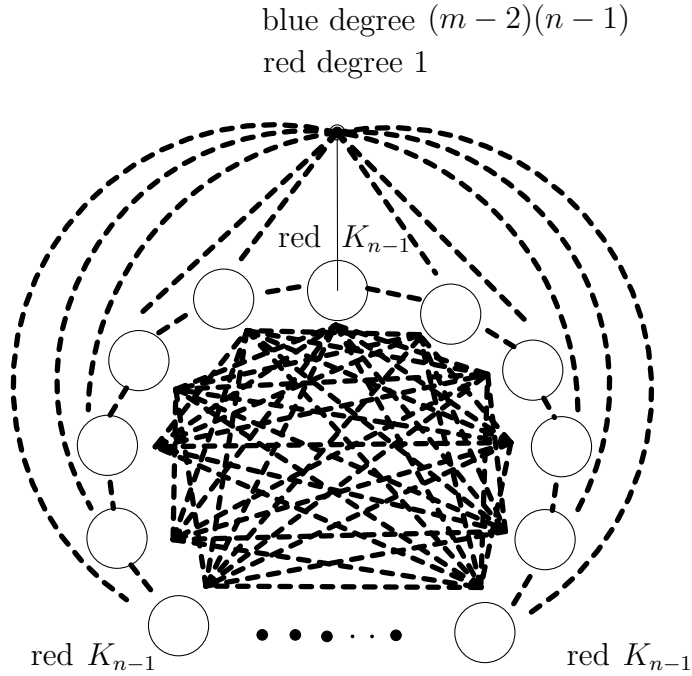
Figure 3. A red $C_n$- free coloring of $K_{(m-1)(n-1)+1}- K_{1,n-2}$ with no blue $K_m$.

\vspace{16pt}

\noindent Hence,  $K_{(m-1)(n-1)+1} - K_{1,n-2} \not \rightarrow (C_n,K_m)$. Therefore, $r_*(C_n, K_m) \geq (m-2)(n-1)+2$.

\vspace{8pt}

\noindent Next to show that, $r_*(C_n, K_5)\leq (m-2)(n-1)+2$, assume that there exists a red $C_n$- free  red/blue coloring of a graph $G=K_{(m-1)(n-1)+1} - K_{1,n-3}$ that contains no blue $K_m$. Let $H$ be the graph obtained by deleting the vertex of degree $(m-2)(n-1)+2$ (say $v$) from $G$ (i.e., $H=G \setminus v$).

\vspace{6pt}

\noindent Then, $H$ is a graph on $(m-1)(n-1)$ vertices such that it contains no red $C_n$ or a blue $K_m$. Therefore, by  lemma \ref{l4} we get that $H$ contains a red $(m-1)K_{n-1}$.  Let $V_1$, $V_2$,...,$V_{m-1}$ denote the  sets of vertices of the $m-1$ connected components of size $n-1$.  In order to avoid a red $C_n$, $v$ can be adjacent to at most one red neighbor in each of the $m-1$ sets $V_1$, $V_2$,...,$V_{m-1}$. We exercise the prerogative of assuming that $V_{m-1}$ represents the connected component that has the least number of vertices connected to $v$. Then as $V_{m-1}$ has at least 2 vertices adjacent to $v$, we can select $v_{m-1} \in V_{m-1}$ such that it is adjacent to $v$ in blue (since if $v$ is adjacent to two vertices of $V_{m-1}$ it will result in a red $C_n$). For each $1 \leq i \leq (m-2)$, let $S_i=\bigcup_{j=1}^{m-2} V_j \setminus V_i$. Clearly, each vertex of $V_i$ ( $1 \leq i \leq m-2$) has at most $m-1$ vertices that have at least one red neighbor in $G[S_i  \cup \{v,v_{m-1}\}]$. As $n>m$, we can conclude that each $V_i$ ( $1 \leq i \leq m-2$)  will have least one vertex  having no red neighbors in $G[S_i  \cup \{v,v_{m-1}\}]$. In other words, we can select $v_i \in V_i$ ($1 \leq i \leq m-2$) lying in the blue neighborhoods of $v_{m-1}$ and $v$  such that $\{v_1,v_2,...,v_{m-2}\}$ induces a blue $K_{m-2}$.  Therefore,  $\{v_1,v_2,v_3,...,v_{m-1},v\}$ will induce a blue $K_m$, a contradiction. Hence the result.

\end{document}